\documentclass[12pt]{amsart}
\usepackage{amssymb,latexsym}
\usepackage{enumerate}

\makeatletter
\@namedef{subjclassname@2010}{%
  \textup{2010} Mathematics Subject Classification}
\makeatother

\newtheorem{theorem}{Theorem}[section]

\newtheorem{corollary}[theorem]{Corollary}

\theoremstyle{definition}

\theoremstyle{remark}
\newtheorem{remark}[theorem]{Remark}
\numberwithin{equation}{section}

\begin{document}

\baselineskip=17pt

\title[Parametrized Borsuk-Ulam problem]{Parametrized Borsuk-Ulam problem for projective space bundles}

\author{Mahender Singh}

\address{Institute of Mathematical Sciences\\ C I T Campus\\ Taramani\\ Chennai 600113\\ India}

\email{mahender@imsc.res.in, mahen51@gmail.com}

\subjclass[2010]{Primary 55M20; Secondary 55R91, 55R25}

\keywords{Characteristic polynomial of a bundle, cohomological dimension, continuity of \v{C}ech cohomology, equivariant map, free involution.}

\begin{abstract}
Let $\pi: E \to B$ be a fiber bundle with fiber having the mod 2 cohomology algebra of a real or a complex projective space and let $\pi^{'}: E^{'} \to B$ be vector bundle such that $\mathbb{Z}_2$ acts fiber preserving and freely on $E$ and $E^{'}-0$, where $0$ stands for the zero section of the bundle $\pi^{'}:E^{'} \to B$. For a fiber preserving $\mathbb{Z}_2$-equivariant map $f:E \to E^{'}$, we estimate the cohomological dimension of the zero set $Z_f = \{x \in E ~|~ f(x)= 0\}.$ As an application, we also estimate the cohomological dimension of the $\mathbb{Z}_2$-coincidence set $A_f=\{x \in E~|~ f(x) = f(T(x)) \}$ of a fiber preserving map $f:E \to E^{'}$.
\end{abstract}

\maketitle

\section{Introduction}
The unit $n$-sphere $\mathbb{S}^n$ is equipped with the antipodal involution given by $x \mapsto -x$. The well known Borsuk-Ulam theorem states that: If $n \geq k$, then for every continuous map $f: \mathbb{S}^n \to \mathbb{R}^k$ there exists a point $x$ in $\mathbb{S}^n$ such that $f(x)= f(-x)$. Over the years there have been several generalizations of the theorem in many directions. We refer the reader to the article \cite{Steinlein} by Steinlein which lists 457 publications concerned with various generalizations of the Borsuk-Ulam theorem.

Jaworowski \cite{Jawo1}, Dold \cite{Dold}, Nakaoka \cite{Nakaoka} and others extended this theorem to the setting of fiber bundles, by considering fiber preserving maps $f : SE \to E^{'}$, where $SE$ denotes the total space of the sphere bundle $SE \to B$ associated to a vector bundle $E \to B$, and $E^{'} \to B$ is other vector bundle. Thus, they parametrized the Borsuk-Ulam theorem, whose general formulation is as follows:
\vspace{1mm}

\textit{Let $G$ be a compact Lie group. Consider a fiber bundle $\pi: E \to B$ and a vector bundle $\pi^{'}: E^{'} \to B$ such that $G$ acts fiber preserving and freely on $E$ and $E^{'}-0$, where $0$ stands for the zero section of the bundle $\pi^{'}:E^{'} \to B$. For a fiber preserving $G$-equivariant map $f:E  \to E^{'}$, the parametrized version of the Borsuk-Ulam theorem deals in estimating the cohomological dimension of the zero set}
$$ Z_f = \{x \in E ~|~ f(x)= 0\}.$$

Such results appeared first in the papers of Jaworowski \cite{Jawo1}, Dold \cite{Dold} and Nakaoka \cite{Nakaoka}. Dold \cite{Dold} and Nakaoka \cite{Nakaoka} defined certain polynomials, which they called the characteristic polynomials, for vector bundles with free $G$-actions ($G = \mathbb{Z}_p$ or $\mathbb{S}^1$) and used them for obtaining such results. The characteristic polynomials were used by Koikara and Mukerjee \cite{Koikara} to prove a parametrized version of the Borsuk-Ulam theorem for bundles whose fiber is a product of spheres, with the free involution given by the product of the antipodal involutions. Recently, Mattos and Santos \cite{Mattos} also used the same technique to obtain parametrized Borsuk-Ulam theorems for bundles whose fiber has the mod $p$ cohomology algebra (with $p>2$) of a product of two spheres with any free $\mathbb{Z}_p$-action and for bundles whose fiber has the rational cohomology algebra of a product of two spheres with any free $\mathbb{S}^1$-action. Jaworowski obtained parametrized Borsuk-Ulam theorems for lens space bundles in \cite{Jawo4} and parametrized Borsuk-Ulam theorems for sphere bundles in \cite{Jawo1, Jawo2, Jawo3}.

The purpose of this paper is to prove parametrized Borsuk-Ulam theorems for bundles whose fiber has the mod 2 cohomology algebra of a real or a complex projective space with any free involution. The theorems are stated in section 4 and proved in section 6. As an application, in section 7, the cohomological dimension of the $\mathbb{Z}_2$-coincidence set of a fiber preserving map is also estimated.

\section{Preliminaries}
Here we recall some basic notions that will be used in later sections. All spaces under consideration will be paracompact Hausdorff spaces and the cohomology used will be the \v{C}ech cohomology with $\mathbb{Z}_2$ coefficients. We will exploit the continuity property of the \v{C}ech cohomology theory, for the details of which we refer to Eilenberg-Steenrod \cite[Chapter X]{Eilenberg}.

We recall that a finitistic space is a paracompact Hausdorff space whose every open covering has a finite dimensional open refinement, where the dimension of a covering is one less than the maximum number of members of the covering which intersect non-trivially (the notion was introduced by Swan in \cite{Swan}). It is a large class of spaces including all compact Hausdorff spaces and all paracompact spaces of finite covering dimension.

For a space $X$, $cohom.dim(X)$ will mean the cohomological dimension of $X$ with respect to $\mathbb{Z}_2$. For basic results of dimension theory, we refer to Nagami \cite{Nagami}. If $G$ is a compact Lie group acting freely on a paracompact Hausdorff space $X$, then $$X \to X/G$$ is a principal $G$-bundle and we can take a classifying map $$X/G  \to B_G$$ for the principal $G$-bundle $X \to X/G$, where $B_G$ is the classifying space of the group $G$. Recall that for $G= \mathbb{Z}_2$, we have $$H^*(B_G; \mathbb{Z}_2)\cong \mathbb{Z}_2[s],$$ where $s$ is a homogeneous element of degree one. We will also use some elementary notions about vector bundles for the details of which we refer to Husemoller \cite{Husemoller}.
\bigskip

\section{Free involutions on projective spaces and their orbit spaces}
We note that odd dimensional real projective spaces admit free involutions. Let $n=2m-1$ with $m \geq 1$. Recall that $\mathbb{R}P^n$ is the orbit space of the antipodal involution on $\mathbb{S}^n$ given by $$(x_1,x_2,..., x_{2m-1},x_{2m}) \mapsto (-x_1,-x_2,..., -x_{2m-1},-x_{2m}).$$ If we denote an element of $\mathbb{R}P^n$ by $[x_1,x_2,..., x_{2m-1},x_{2m}]$, then the map $\mathbb{R}P^n \to \mathbb{R}P^n$ given by $$[x_1,x_2,..., x_{2m-1},x_{2m}] \mapsto [-x_2,x_1,..., -x_{2m},x_{2m-1}]$$ defines an involution. If $$[x_1,x_2,..., x_{2m-1},x_{2m}]=[-x_2,x_1,..., -x_{2m},x_{2m-1}],$$ then $$(-x_1,-x_2,..., -x_{2m-1},-x_{2m})=(-x_2,x_1,..., -x_{2m},x_{2m-1}),$$ which gives $x_1=x_2=... =x_{2m-1}=x_{2m}=0$, a contradiction. Hence, the involution is free.

Similarly, the complex projective space $\mathbb{C}P^n$ admit free involutions when $n \geq 1$ is odd. Recall that $\mathbb{C}P^n$ is the orbit space of the free $\mathbb{S}^1$-action on $\mathbb{S}^{2n+1}$ given by $$(z_1,z_2,..., z_n,z_{n+1}) \mapsto (\zeta z_1,\zeta z_2,..., \zeta z_n,\zeta z_{n+1}) ~ \textrm{for}~ \zeta \in \mathbb{S}^1.$$ If we denote an element of $\mathbb{C}P^n$ by $[z_1,z_2,..., z_n,z_{n+1}]$, then the map $$[z_1,z_2,..., z_n,z_{n+1}] \mapsto [-\overline{z}_2,\overline{z}_1,..., -\overline{z}_{n+1},\overline{z}_n]$$ defines an involution on $\mathbb{C}P^n$. If $$[z_1,z_2,..., z_n,z_{n+1}] = [-\overline{z}_2,\overline{z}_1,..., -\overline{z}_{n+1},\overline{z}_n],$$ then $$(\lambda z_1,\lambda z_2,..., \lambda z_n,\lambda z_{n+1}) = (-\overline{z}_2,\overline{z}_1,..., -\overline{z}_{n+1},\overline{z}_n)$$ for some $\lambda \in \mathbb{S}^1$, which gives $z_1=z_2=... =z_n=z_{n+1}=0$, a contradiction. Hence, the involution is free.

We write $X \simeq_2 \mathbb{R}P^n$ if $X$ is a space having the mod 2 cohomology algebra of $\mathbb{R}P^n$. Similarly, we write $X \simeq_2 \mathbb{C}P^n$ if $X$ is a  space having the mod 2 cohomology algebra of $\mathbb{C}P^n$.

Recently, Singh \cite{hksingh} determined completely the mod 2 cohomology algebra of orbit spaces of free involutions on mod 2 cohomology real and complex projective spaces. Using the Leray spectral sequence associated to the Borel fibration $$X\hookrightarrow X_G \longrightarrow B_G,$$ they proved the following results.

\begin{theorem}
If $G=\mathbb{Z}_2$ acts freely on a finitistic space $X \simeq_2 \mathbb{R}P^n$, where $n$ is odd, then $$H^*(X/G; \mathbb{Z}_2)\cong \mathbb{Z}_2[u,v]/\langle u^2, v^{\frac{n+1}{2}}\rangle,$$
where deg($u$)=$1$ and deg($v$)=$2$.
\end{theorem}

\begin{theorem}
If $G=\mathbb{Z}_2$ acts freely on a finitistic space $X \simeq_2 \mathbb{C}P^n$, where $n$ is odd, then $$H^*(X/G; \mathbb{Z}_2) \cong \mathbb{Z}_2[u,v]/\langle u^3, v^{\frac{n+1}{2}}\rangle,$$
where deg($u$)=$1$ and deg($v$)=$4$.
\end{theorem}

\begin{remark}
It is easy to see that, when $n$ is even, then $\mathbb{Z}_2$ cannot act freely on a finitistic space $X \simeq_2 \mathbb{R}P^n$ or $\mathbb{C}P^n$. For, if $n$ is even, then the Euler characteristic
\begin{displaymath}
\chi(X) = \left\{ \begin{array}{ll}
1 & \textrm{when $X \simeq_2 \mathbb{R}P^n$}\\
n+1 & \textrm{when $X \simeq_2 \mathbb{C}P^n$.}\\
\end{array} \right.
\end{displaymath}
But for a free involution, $\chi(X^{\mathbb{Z}_2})=0$ and hence the Floyd's Euler characteristic formula \cite[p.145]{Bredon}
$$\chi(X)+ \chi(X^{\mathbb{Z}_2})= 2 \chi(X/{\mathbb{Z}_2})$$
gives a contradiction.
\end{remark}

\begin{remark}
Let $X \simeq_2 \mathbb{H}P^n$ be a finitistic space, where $\mathbb{H}P^n$ is the quaternionic projective space. For $n=1$, $X \simeq_2 \mathbb{S}^4$, which is dealt in \cite{Dold}. For $n \geq 2$, there is no free involution on $X$, which follows from the stronger fact that such spaces have the fixed point property.
\end{remark}

\begin{remark}
Let $X \simeq_2 \mathbb{O}P^2$ be a finitistic space, where $\mathbb{O}P^2$ is the Cayley projective plane. Note that $H^*(\mathbb{O}P^2; \mathbb{Z}_2) \cong \mathbb{Z}_2[u]/\langle u^3 \rangle$, where $u$ is a homogeneous element of degree 8. Just as in Remark 3.3, it follows from the Floyd's Euler characteristic formula that there is no free involution on $X$.
\end{remark}
\bigskip

\section{Statements of theorems}
Let $X \simeq_2 \mathbb{R}P^n$ be a finitistic space. Let $(X, E, \pi, B)$ be a fiber bundle with a fiber preserving free $\mathbb{Z}_2$-action such that the quotient bundle $(X/G, \overline{E}, \overline{\pi}, B)$ has a cohomology extension of the fiber, that is, there is a $\mathbb{Z}_2$-module homomorphism of degree zero
$$\theta: H^*(X/G; \mathbb{Z}_2) \to H^*(\overline{E}; \mathbb{Z}_2)$$
such that for any $b \in B$, the composition
$$H^*(X/G; \mathbb{Z}_2) \stackrel{\theta}\rightarrow H^*(\overline{E}; \mathbb{Z}_2) \stackrel{i_b^*}\rightarrow H^*((X/G)_b; \mathbb{Z}_2)$$ is an isomorphism, where $$i_b:(X/G)_b \hookrightarrow \overline{E}$$ is the inclusion of the fiber over $b$ (see \cite[p.372]{Bredon}). This condition on the bundle is assumed so that the Leray-Hirsch theorem can be applied (see \cite[Chapter VII, Theorem 1.4]{Bredon}). Now consider a $k$-dimensional vector bundle $$\pi^{'}: E^{'}  \to B$$ with a fiber preserving $\mathbb{Z}_2$-action on $E^{'}$ which is free on $E^{'}-0$. Let $$f: E \to E^{'}$$ be a fiber preserving $\mathbb{Z}_2$-equivariant map. Define $$Z_f = \{x \in E ~|~ f(x)= 0\}$$ and $$\overline{Z_{f}}= Z_f/{\mathbb{Z}_2},$$ the quotient by the free $\mathbb{Z}_2$-action induced on $Z_f$.

Let $H^*(B)[x,y]$ be the polynomial ring over $H^*(B)$ in the indeterminates $x$ and $y$. For the bundle $(X \simeq_2 \mathbb{R}P^n, E, \pi, B)$, in section 5, we will define the characteristic polynomials $W_1(x,y)$ and $W_2(x,y)$ in $H^*(B)[x,y]$ and we will show that $H^*(\overline{E})$ and $H^*(B)[x,y]/\langle W_1(x,y), W_2(x,y)\rangle$ are isomorphic as $H^*(B)$-modules. Therefore, each polynomial $q(x,y)$ in $H^*(B)[x,y]$ defines an element of $H^*(\overline{E})$, which we will denote by $q(x, y)|_{\overline{E}}$. We will denote by $$q(x, y)|_{\overline{Z_f}}$$ the image of $q(x,y)|_{\overline{E}}$ by the $H^*(B)$-homomorphism $$i^*: H^*(\overline{E}) \to H^*(\overline{Z_f}),$$ where $i^*$ is the map induced by the inclusion $i: \overline{Z_f} \hookrightarrow \overline{E}$. In a similar way, we will define the characteristic polynomial $W^{'}(x)$ for the vector bundle $\pi^{'}: E^{'} \to B$. With the above hypothesis and notations, we prove the following results for the real case.

\begin{theorem}
Let $X \simeq_2 \mathbb{R}P^n$ be a finitistic space. If $q(x,y)$ in $H^*(B)[x,y]$ is a polynomial such that $q(x,y)|_{\overline{Z_f}}=0$, then there are polynomials $r_1(x,y)$ and $r_2(x,y)$ in $H^*(B)[x,y]$ such that $$q(x,y)W^{'}(x) = r_1(x,y)W_1(x,y)+ r_2(x,y)W_2(x,y)$$ in the ring $H^*(B)[x,y]$, where $W^{'}(x)$, $W_1(x,y)$ and $W_2(x,y)$ are the characteristic polynomials.
\end{theorem}

As a corollary, just as in \cite{Dold}, we have the following parametrized version of the Borsuk-Ulam theorem.

\begin{corollary}
Let $X \simeq_2 \mathbb{R}P^n$ be a finitistic space. If the fiber dimension of $E^{'} \to B$ is $k$, then $q(x,y)|_{\overline{Z_f}}\neq 0$ for all non zero polynomials $q(x,y)$ in $H^*(B)[x,y]$, whose degree in $x$ and $y$ is less than $(n-k+1)$. Equivalently, the $H^*(B)$-homomorphism
$$\bigoplus_{i+j=0}^{n-k}H^*(B)x^iy^j \to H^*(\overline{Z_f})$$
given by $x^i \to x^i|_{\overline{Z_f}}$ and $y^j \to y^j|_{\overline{Z_f}}$ is a monomorphism. As a result, if $n \geq k$, then 
$$cohom.dim (Z_f) \geq cohom.dim (B) + (n-k).$$
\end{corollary}
\bigskip

Let $X \simeq_2 \mathbb{C}P^n$ be a finitistic space. Just as in the real case, we will define the characteristic polynomials $W_1(x,y)$ and $W_2(x)$ for the bundle $(X \simeq_2 \mathbb{C}P^n, E, \pi, B)$ and show that $H^*(\overline{E})$ and $H^*(B)[x,y]/\langle W_1(x,y), W_2(x)\rangle$ are isomorphic as $H^*(B)$-modules. With similar hypothesis and notations as in the real case, we prove the following results for the complex case.

\begin{theorem}
Let $X \simeq_2 \mathbb{C}P^n$ be a finitistic space. If $q(x,y)$ in $H^*(B)[x,y]$ is a polynomial such that $q(x, y)|_{\overline{Z_f}}=0$, then there are polynomials $r_1(x,y)$ and $r_2(x,y)$ in $H^*(B)[x,y]$ such that $$q(x,y)W^{'}(x) = r_1(x,y)W_1(x,y)+ r_2(x,y)W_2(x)$$ in the ring $H^*(B)[x,y]$, where $W^{'}(x)$, $W_1(x,y)$ and $W_2(x)$ are the characteristic polynomials.
\end{theorem}

\begin{corollary}
Let $X \simeq_2 \mathbb{C}P^n$ be a finitistic space. If the fiber dimension of $E^{'} \to B$ is $k$, then $q(x,y)|_{\overline{Z_f}} \neq 0$ for all non zero polynomials $q(x,y)$ in $H^*(B)[x,y]$, whose degree in $x$ and $y$ is less than $(2n-k+2)$. Equivalently, the $H^*(B)$-homomorphism 
$$\bigoplus_{i+j=0}^{2n-k+1}H^*(B)x^iy^j \to H^*(\overline{Z_f})$$
given by $x^i \to x^i|_{\overline{Z_f}}$ and $y^j \to y^j|_{\overline{Z_f}}$ is a monomorphism. As a result, if $2n \geq k$, then
$$cohom.dim (Z_f) \geq cohom.dim (B) + (2n-k+1).$$
\end{corollary}
\bigskip

\section{Characteristic polynomials for bundles}
Let $(X, E, \pi, B)$ be a fiber bundle with a fiber preserving free $\mathbb{Z}_2$-action such that the quotient bundle $(X/G, \overline{E}, \overline{\pi}, B)$ has a cohomology extension of the fiber. With this hypothesis, we now proceed to define the characteristic polynomials for the bundles. We deal the real and the complex case separately.
\bigskip

\noindent \textbf{When $X \simeq_2 \mathbb{R}P^n$.}\\
Let $G=\mathbb{Z}_2$ act freely on a finitistic space $X \simeq_2 \mathbb{R}P^n$. Then, $n$ is odd and by theorem 3.1, $H^*(X/G;\mathbb{Z}_2)$ is a free graded algebra generated by the elements
$$1,~ u,~ v,~ uv,..., v^{\frac{n-1}{2}},~ uv^{\frac{n-1}{2}},$$
subject to the relations $u^2= 0$ and $v^{\frac{n+1}{2}}=0$, where $u \in H^1(X/G;\mathbb{Z}_2)$ and  $v \in H^2(X/G;\mathbb{Z}_2)$. Let $(X \simeq_2 \mathbb{R}P^n, E, \pi, B)$ be a bundle with the hypothesis of section 4. By the Leray-Hirsch theorem, there exist elements $a \in H^1(\overline{E})$ and $b \in H^2(\overline{E})$ such that the restriction to a typical fiber $$j^*:H^*(\overline{E}) \to H^*(X/G)$$ maps $a \mapsto u$ and $b \mapsto v$. Note that $H^*(\overline{E})$ is a
$H^*(B)$-module, via the induced homomorphism $\overline{\pi}^*$ and is generated by the basis
\begin{equation}\label{eqn1}
1,~ a,~ b,~ ab,..., b^{\frac{n-1}{2}},~ ab^{\frac{n-1}{2}}.
\end{equation}
We can express the element $b^{\frac{n+1}{2}} \in H^{n+1}(\overline{E})$ in terms of the basis \eqref{eqn1}. Therefore, there exist unique elements $w_i \in H^i(B)$ such that
$$ b^{\frac{n+1}{2}} = w_{n+1} + w_na + w_{n-1}b + \cdots  + w_2b^{\frac{n-1}{2}} + w_1ab^{\frac{n-1}{2}}.$$
Similarly, we express the element $a^2 \in H^2(\overline{E})$ as $$a^2= \nu_2 + \nu_1a +\alpha b,$$ where $\nu_i \in H^i(B)$ and $ \alpha \in \mathbb{Z}_2$ are unique elements.
The characteristic polynomials in the indeterminates $x$ and $y$, of degrees respectively 1 and 2, associated to the fiber bundle $(X \simeq_2 \mathbb{R}P^n, E, \pi, B)$ are defined by
$$ W_1(x,y)= w_{n+1} + w_nx + w_{n-1}y + \cdots  + w_2y^{\frac{n-1}{2}} + w_1xy^{\frac{n-1}{2}} + y^{\frac{n+1}{2}}$$
$$\textrm{and}~~W_2(x,y)= \nu_2 + \nu_1x +\alpha y+ x^2.$$
On substituting the values for the indeterminates $x$ and $y$, we obtain the following homomorphism of $H^*(B)$-algebras
$$\sigma: H^*(B)[x,y] \to H^*(\overline{E})$$
given by $(x,y) \mapsto (a,b)$. The $Ker(\sigma)$ is the ideal generated by the polynomials $W_1(x,y)$ and $W_2(x,y)$ and hence
\begin{equation}\label{eqn2}
H^*(B)[x,y]/\langle W_1(x,y), W_2(x,y) \rangle \cong H^*(\overline{E}).
\end{equation}
\bigskip

\noindent \textbf{When $X \simeq_2 \mathbb{C}P^n$.}\\
Since this case is similar, we present it rather briefly. Let $G=\mathbb{Z}_2$ act freely on a finitistic space $X \simeq_2 \mathbb{C}P^n$. Then, $n$ is odd and by theorem 3.2, $H^*(X/G;\mathbb{Z}_2)$ is a free graded algebra generated by the elements
$$1,~ u,~ u^2,~ v,~ uv,..., v^{\frac{n-1}{2}},~ uv^{\frac{n-1}{2}},~ u^2v^{\frac{n-1}{2}},$$
subject to the relations $u^3= 0$ and $v^{\frac{n+1}{2}}=0$, where $u \in H^1(X/G;\mathbb{Z}_2)$ and $v \in H^4(X/G;\mathbb{Z}_2)$. By the Leray-Hirsch theorem, there exist elements $a \in H^1(\overline{E})$ and $b \in H^4(\overline{E})$ such that the restriction to a typical fiber $$j^*:H^*(\overline{E}) \to H^*(X/G)$$ maps $a \mapsto u$ and $b \mapsto v$. Note that $H^*(\overline{E})$ is a $H^*(B)$-module and is generated by the basis
\begin{equation}\label{eqn3}
1,~ a,~ a^2,~ b,~ ab,..., b^{\frac{n-1}{2}},~ ab^{\frac{n-1}{2}},~ a^2b^{\frac{n-1}{2}}.
\end{equation}
We write $b^{\frac{n+1}{2}} \in H^{2n+2}(\overline{E})$ in terms of the basis \eqref{eqn3}. Thus, there exist unique elements $w_i \in H^i(B)$ such that
$$ b^{\frac{n+1}{2}} = w_{2n+2} + w_{2n+1}a + w_{2n}a^2 + \cdots  + w_2a^2b^{\frac{n-1}{2}}.$$
Similarly, we write the element $a^3 \in H^3(\overline{E})$ as $$a^3 = \nu_3 + \nu_2a + \nu_1a^2,$$ where $\nu_i \in H^i(B)$ are unique elements.
The characteristic polynomials in the indeterminates $x$ and $y$, of degrees respectively 1 and 4, associated to the fiber bundle $(X \simeq_2 \mathbb{C}P^n, E, \pi, B)$ are defined by
$$ W_1(x,y)= w_{2n+2} + w_{2n+1}x + w_{2n}x^2 + \cdots  + w_2x^2y^{\frac{n-1}{2}} + y^{\frac{n+1}{2}}$$
$$\textrm{and}~~W_2(x)= \nu_3 + \nu_2x + \nu_1x^2 + x^3.$$
This gives a homomorphism of $H^*(B)$-algebras
$$\sigma: H^*(B)[x,y] \to H^*(\overline{E})$$
given by $(x,y) \mapsto (a,b)$ and having $Ker(\sigma)$ as the ideal generated by the polynomials
$W_1(x,y)$ and $W_2(x)$. Hence
\begin{equation}\label{eqn4}
H^*(B)[x,y]/\langle W_1(x,y), W_2(x) \rangle \cong H^*(\overline{E}).
\end{equation}
\bigskip

\noindent \textbf{Characteristic polynomial for the bundle $\pi{'}: E^{'}  \to B$.}\\
Now we define the characteristic polynomial associated to the $k$-dimensional vector bundle $\pi{'}: E^{'}  \to B$ with fiber preserving $\mathbb{Z}_2$-action on $E^{'}$ which is free on $E^{'}-0$. Let $SE^{'}$ denote the total
space of sphere bundle of $\pi{'}: E^{'}  \to B$ . Since the action is free on $SE^{'}$, we obtain the projective space bundle $(\mathbb{R}P^{k-1}, \overline{SE^{'}}, \overline{\pi^{'}}, B)$ and the principal $\mathbb{Z}_2$-bundle $SE^{'} \to \overline{SE^{'}}$.
We know that $$H^*(\mathbb{R}P^{k-1}; \mathbb{Z}_2) \cong \mathbb{Z}_2[u^{'}]/ \langle {u^{'}}^k\rangle,$$ where $u^{'}= g^*(s)$, $s \in H^1(B_G)$ and $g: \mathbb{R}P^{k-1} \to B_G$ is a classifying map for the principal $\mathbb{Z}_2$-bundle $\mathbb{S}^{k-1} \to \mathbb{R}P^{k-1}$. Let $h: \overline{SE^{'}} \to  B_G$ be a classifying map for the principal $\mathbb{Z}_2$-bundle $SE^{'} \to \overline{SE^{'}}$ and let $a^{'}= h^*(s) \in H^1(\overline{SE^{'}})$. Now the $\mathbb{Z}_2$-module homomorphism $$\theta^{'}:H^*(\mathbb{R}P^{k-1}) \to H^*(\overline{SE^{'}})$$ given by $u^{'} \mapsto a^{'}$ is a cohomology extension of the fiber. Again, by the Leray-Hirsch theorem $H^*(\overline{SE^{'}})$ is a $H^*(B)$-module via the induced homomorphism ${\overline{\pi^{'}}}^*$ and is generated by the basis
$$1, a^{'}, {a^{'}}^2, ..., {a^{'}}^{k-1}.$$
We write ${a^{'}}^{k}\in H^k(\overline{SE^{'}})$ as 
$${a^{'}}^{k} = w_k^{'} + w_{k-1}^{'}a^{'}+ \cdots + w_1^{'}{a^{'}}^{k-1},$$
where $w_i^{'} \in H^i(B)$ are unique elements. Now the characteristic polynomial in the indeterminate $x$ of degree 1, associated to the vector bundle $\pi{'}: E^{'}  \to B$ is defined as
$$ W^{'}(x)= w_k^{'} + w_{k-1}^{'}x + \cdots + w_1^{'}x^{k-1} +x^k.$$
By similar arguements as used above, we have the following isomorphism of $H^*(B)$-algebras
$$ H^*(B)[x]/\langle W^{'}(x) \rangle \cong H^*(\overline{SE^{'}})$$
given by $x \mapsto a^{'}.$
\bigskip

\section{Proofs of theorems}
We first prove our results for the real case.
\bigskip

\noindent \textit{Proof of Theorem 4.1.}
Let $q(x,y)$ in $H^*(B)[x,y]$ be a polynomial such that $q(x,y)|_{\overline{Z_f}} = 0$. It follows from the continuity property of the \v{C}ech cohomology theory, that there is an open subset $V \subset \overline{E}$ such that $\overline{Z_f} \subset V$ and $q(x,y)|_V = 0$. Consider the long exact cohomology sequence for the pair $(\overline{E}, V)$, namely,
$$\cdots \to H^*(\overline{E}, V) \stackrel{j_1^*}{\rightarrow} H^*(\overline{E}) \rightarrow H^*(V) \rightarrow H^*(\overline{E}, V) \rightarrow \cdots.$$
By exactness, there exists $ \mu \in H^*(\overline{E}, V)$ such that $j_1^*(\mu)= q(x,y)|_{\overline{E}}$, where $$j_1: \overline{E} \to (\overline{E}, V)$$ is the natural inclusion. The $\mathbb{Z}_2$-equivariant map $$f: E \to E^{'}$$ gives the map $$\overline{f}:\overline{E}-\overline{Z_f} \rightarrow \overline{E^{'}}-0.$$ The induced map $${\overline{f}}^*: H^*(\overline{E^{'}}-0) \to H^*(\overline{E}-\overline{Z_f})$$ is a $H^*(B)$-homomorphism. Also we have $W^{'}(a^{'})=0$. Therefore,
$$W^{'}(x)|_{\overline{E}-\overline{Z_f}}= W^{'}(a)= W^{'}\big( {\overline{f}}^* (a^{'})\big)= {\overline{f}}^* \big( W^{'}(a^{'}) \big)=0.$$
Now consider the long exact cohomology sequence for the pair $(\overline{E}, \overline{E}-\overline{Z_f})$, that is,
$$\cdots \to H^*(\overline{E}, \overline{E}-\overline{Z_f}) \stackrel{j_2^*}{\rightarrow} H^*(\overline{E}) \rightarrow H^*(\overline{E}-\overline{Z_f}) \rightarrow H^*(\overline{E}, \overline{E}-\overline{Z_f}) \rightarrow \cdots.$$
By exactness, there exists $ \lambda \in H^*(\overline{E}, \overline{E}-\overline{Z_f})$ such that $j_2^*(\lambda)= W^{'}(x)|_{\overline{E}}$, where $$j_2: \overline{E} \to (\overline{E}, \overline{E}-\overline{Z_f})$$ is the natural inclusion. Thus,
$$q(x,y)W^{'}(x)|_{\overline{E}}= j_1^*(\mu)\smile j_2^*(\lambda)=j^*(\mu \smile \lambda)$$
by the naturality of the cup product. But, $$\mu \smile \lambda \in H^*\big(\overline{E}, V \cup (\overline{E}-\overline{Z_f})\big)=H^*(\overline{E}, \overline{E})=0$$ and hence $q(x,y)W^{'}(x)|_{\overline{E}}=0$. Therefore, by equation \eqref{eqn2}, there exist polynomials $r_1(x,y)$ and $r_2(x,y)$ in $H^*(B)[x,y]$ such that $$q(x,y)W^{'}(x) = r_1(x,y)W_1(x,y)+ r_2(x,y)W_2(x,y)$$ in the ring $H^*(B)[x,y]$. This proves the theorem. \hfill $\Box$
\bigskip

\noindent \textit{Proof of Corollary 4.2.}
Let $q(x,y)$ in $H^*(B)[x,y]$ be a non zero polynomial such that deg$(q(x,y)) < (n-k+1)$. If $q(x,y)|_{\overline{Z_f}}=0$, then by theorem 4.1, we have $$q(x,y)W^{'}(x) = r_1(x,y)W_1(x,y) + r_2(x,y)W_2(x,y)$$ in the ring $H^*(B)[x,y]$ for some polynomials $r_1(x,y)$ and $r_2(x,y)$ in $H^*(B)[x,y]$. Note that deg$(W^{'}(x))=k$, deg$(W_1(x,y))=n+1$ and deg$(W_2(x,y))=2$. Since
$$\textrm{deg}(q(x,y)) + k = \textrm{max} \{ \textrm{deg}(r_1(x,y)) + n+1, ~\textrm{deg}(r_2(x,y)) + 2 \},$$
we have $$\textrm{deg}(q(x,y)) + k \geq  \textrm{deg}(r_1(x,y)) + n+1.$$
Taking deg$(r_1(x,y)) =0$, this gives deg$(q(x,y)) + k \geq n+1$ and hence deg$(q(x,y)) \geq (n-k+1)$, which is a contradiction. Hence $q(x,y)|_{\overline{Z_f}}\neq 0$.
Equivalently, the $H^*(B)$-homomorphism 
$$\bigoplus_{i+j=0}^{n-k}H^*(B)x^iy^j \to H^*(\overline{Z_f})$$
given by $x^i \to x^i|_{\overline{Z_f}}$ and $y^j \to y^j|_{\overline{Z_f}}$ is a monomorphism. As a result, if $n \geq k$, then
$$cohom.dim (Z_f) \geq cohom.dim (B) + (n-k),$$
since $cohom.dim (Z_f) \geq cohom.dim (\overline{Z_f})$ by \cite[Proposition A.11]{Quillen}. \hfill $\Box$

\begin{remark}
If $B$ is a point in the above corollary, then for any $\mathbb{Z}_2$-equivariant map $$f:X \simeq_2 \mathbb{R}P^n \to \mathbb{R}^k,$$ where $n \geq k$, we have $cohom.dim (Z_f) \geq (n-k)$.
\end{remark}

Next we prove our results for the complex case.\\

\noindent \textit{Proof of Theorem 4.3.}
Let $q(x,y)$ in $H^*(B)[x,y]$ be a polynomial such that $q(x,y)|_{\overline{Z_f}} = 0$. By similar arguements as used in the proof of theorem 4.1, we conclude that $q(x,y)W^{'}(x)|_{\overline{E}}=0$. Therefore, by equation \eqref{eqn4}, there exist polynomials $r_1(x,y)$ and $r_2(x,y)$ in $H^*(B)[x,y]$ such that $$q(x,y)W^{'}(x) = r_1(x,y)W_1(x,y)+ r_2(x,y)W_2(x)$$ in the ring $H^*(B)[x,y]$. This proves the theorem.\hfill $\Box$
\bigskip

\noindent \textit{Proof of Corollary 4.4.}
Let $q(x,y)$ in $H^*(B)[x,y]$ be a non zero polynomial such that deg$(q(x,y)) < (2n-k+2)$. If $q(x,y)|_{\overline{Z_f}}=0$, then by theorem 4.3, we have $$q(x,y)W^{'}(x) = r_1(x,y)W_1(x,y) + r_2(x,y)W_2(x)$$ in the ring $H^*(B)[x,y]$ for some polynomials $r_1(x,y)$ and $r_2(x,y)$ in $H^*(B)[x,y]$. Note that deg$(W^{'}(x))=k$, deg$(W_1(x,y))=2n+2$ and deg$(W_2(x))=3$. Since
$$\textrm{deg}(q(x,y)) + k = \textrm{max} \{ \textrm{deg}(r_1(x,y)) + 2n+2, ~\textrm{deg}(r_2(x,y)) + 3 \},$$
we have $$\textrm{deg}(q(x,y)) + k \geq  \textrm{deg}(r_1(x,y)) + 2n+2.$$
Taking deg$(r_1(x,y)) =0$, this gives deg$(q(x,y)) + k \geq 2n+2$ and hence deg$(q(x,y)) \geq (2n-k+2)$, which is a contradiction. Hence $q(x,y)|_{\overline{Z_f}}\neq 0$.
Equivalently, the $H^*(B)$-homomorphism
$$\bigoplus_{i+j=0}^{2n-k+1}H^*(B)x^iy^j \to H^*(\overline{Z_f})$$
given by $x^i \to x^i|_{\overline{Z_f}}$ and $y^j \to y^j|_{\overline{Z_f}}$ is a monomorphism. As a result, if $2n \geq k$, then
$$cohom.dim (Z_f) \geq cohom.dim (B) + (2n-k+1). \hspace{20mm} \Box$$

\begin{remark}
If $B$ is a point in the above corollary, then for any $\mathbb{Z}_2$-equivariant map $$f:X \simeq_2 \mathbb{C}P^n \to \mathbb{R}^k,$$ where $2n \geq k$, we have $cohom.dim (Z_f) \geq (2n-k+1)$.
\end{remark}
\bigskip

\section{Application to $\mathbb{Z}_2$-coincidence sets}
Let $(X, E, \pi, B)$ be a fiber bundle with the hypothesis of section 4. Let $E^{''} \to B$ be a $k$-dimensional vector bundle and let $f: E  \to E^{''}$ be a fiber preserving map. Here we do not assume that $E^{''}$ has an involution. Even if $E^{''}$ has an involution, $f$ is not assumed to be $\mathbb{Z}_2$-equivariant. If $T: E  \to E$ is a generator of the $\mathbb{Z}_2$-action, then the $\mathbb{Z}_2$-coincidence set of $f$ is defined as
$$A_f=\{x \in E~|~ f(x) = f(T(x)) \}.$$
Let $V=E^{''} \oplus E^{''}$ be the Whitney sum of two copies of $E^{''} \to B$. Then $\mathbb{Z}_2$ acts on $V$ by permuting the coordinates. This action has the diagonal $D$ in $V$ as the fixed point set. Note that $D$ is a $k$-dimensional sub-bundle of $V$ and the orthogonal complement $D^{\perp}$ of $D$ is also a $k$-dimensional sub-bundle of $V$. Also note that $D^{\perp}$ is $\mathbb{Z}_2$-invariant and has a $\mathbb{Z}_2$-action which is free outside the zero section. Consider the $\mathbb{Z}_2$-equivariant map $f^{'}: E \to V$ given by $$f^{'}(x)= \big(f(x), f(T(x))\big).$$ The linear projection along the diagonal defines a $\mathbb{Z}_2$-equivariant fiber preserving map $g:V \to D^{\perp}$ such that $g(V-D) \subset D^{\perp}-0$, where 0 is the zero section of $D^{\perp}$. Let $h = g \circ f^{'}$ be the composition
$$(E, E-A_f) \to (V, V-D) \to (D^{\perp}, D^{\perp}-0).$$
Note that $$Z_h = h^{-1}(0) = {f^{'}}^{-1}(D)=A_f$$ and $h:E \to D^{\perp}$ is a fiber preserving $\mathbb{Z}_2$-equivariant map.

Applying corollary 4.2 to $h$, we have

\begin{theorem}
If $X \simeq_2 \mathbb{R}P^n$ is a finitistic space, then $$cohom.dim (A_f) \geq cohom.dim(B) + (n-k).$$
\end{theorem}

Similarly, applying corollary 4.4 to $h$, we have
\begin{theorem}
If $X \simeq_2 \mathbb{C}P^n$ is a finitistic space, then $$cohom.dim (A_f) \geq cohom.dim(B) + (2n-k+1).$$
\end{theorem}

\end{document}